\documentclass[11pt]{article}
\usepackage{amsfonts}

\begin{document}

\def\det{\mbox{det}}
\def\Spin{\mbox{Spin}}
\def\Tr{\mbox{Tr \-}}
\def\inj{\mbox{inj \-}}
\def\R{{\mathbb R}}
\def\C{{\mathbb C}}
\def\H{{\mathbb H}}
\def\Ca{{\mathbb Ca}}
\def\Z{{\mathbb Z}}
\def\N{{\mathbb N}}
\def\Q{{\mathbb Q}}
\def\Ad{\mbox{Ad \-}}
\def\k{{\bf k}}
\def\l{{\bf l}}
\def\diag{\mbox{diag \-}}
\def\a{{\alpha}}

\author{Ya.~V.~Bazaikin and E.~G.~Malkovich}

\title{The $Spin(7)$-structures on complex line bundles and explicit
Riemannian metrics with $SU(4)$-holonomy. \footnote{This work was
supported by the RFBR (grant 09-01-00598-a); also first author was
supported by the RFBR (grant 09-01-12130-ofi-m); second author was
supported by the RFBR (grant 09-01-00142-a) and Federal Target Grant
``Scientific and educational personnel of innovation Russia'' for
2009-2013 (government contract No. 02.740.11.0429)}}

\date{}

\maketitle

\begin{abstract}

We completely explore the system of ODE's which is equivalent to the
existence of a parallel $Spin(7)$-structure on the cone over a
7-dimensional 3-Sasakian manifold. The one-dimensional family of
solutions of this system is constructed. The solutions of this
family correspond to metrics with holonomy $SU(4)$ which generalize
the Calabi metrics.


\end{abstract}

\sloppy

\section[]{Introduction.}

The first example of a complete Riemannian metric with holonomy
$SU(n)$ was found explicitly by Calabi \cite{Calabi} in terms of
algebraic functions. The Calabi metric is defined on an appropriate
line complex bundle over arbitrary Einstein-K\"{a}hler manifold $F$.
In the case $F = \C P^{n-1}$ the corresponding Calabi metric is
asymptotically locally euclidean (ALE), in other cases it is
asymptotically conical (AC). In the same paper \cite{Calabi} Calabi
studies hyperk\"{a}hler metrics and explicitly constructs complete
Riemannian metrics with holonomy group $Sp(m)$ on $T^* \C P^m$ ---
the very first explicit examples of hyperk\"{a}hler complete metric.

Here we construct in explicit, algebraic terms one-dimensional
family of complete Riemannian metrics ``connecting'' these two
Calabi's metrics in the space of special K\"{a}hler 8-dimensional
metrics for $F$ to be equal the space of complex $3$-flags in
$\C^3$. In this situation the 4-dimensional quaternion-K\"{a}hler
manifold ${\cal O}$, naturally associated with $F$, admits invariant
``splitting'' of the tangent space. This allows to consider one more
parameter for metric deformation and to solve corresponding system
of ODE's in elementary functions.

The interest in explicit metrics with special holonomies
(particularly, in special K\"{a}hler metrics) is rather high because
there are few examples of such metrics. For instance, Joyce in
\cite[8.2.5]{Joyce} conjectured, that any other (with the exception
of Calabi's metrics for $F=\C P^{n-1}$) ALE-metrics are
``transcendent'', i.e. can not be described in algebraic form. We
should notice that the metrics constructed in this paper are
AC-metrics, but not ALE-metrics, so we didn't make a counterexample
for Joyce's conjecture. As far as we know our metrics bring the
first example of continual family of the Riemannian metrics with
holonomy $SU(n)$, $n \geq 3$, described in elementary functions.

More precisely, let $M=SU(3)/U(1)_{1,1,-2}$ be the Aloff-Wallach
space admitting structure of $3$-Sasakian 7-dimensional manifold. On
$\bar{M}=M \times {\mathbb R}_+$ consider Riemannian metric of the
following form:
$$
dt^2+ A_1(t)^2 \eta_1^2 + A_2(t)^2 \eta_2^2 + A_3(t)^2 \eta_3^2+
B(t)^2 (\eta_4^2+\eta_5^2) + C(t)^2 (\eta_6^2+\eta_7^2),\eqno{(*)}
$$
where $t$ is coordinate in ${\mathbb R}_+$, $\eta_i$ are orthogonal
co-frame on $M$ compatible with $3$-Sasakian structure (details can
be found in the next section). Conic singularity (at $t=0$) of the
$\bar{M}$ can be resolved in the following way: every circle
generated by covector $\eta_1$ collapses into a point when $t
\rightarrow 0$. Manifold obtained after factorizing by ${\mathbb
Z}_2$ is diffeomorphic to $H/ \Z_2$  --- square of the canonical
complex line bundle over space of complex $3$-flags in $\C^3$.

{\bf Theorem 1.} {\it For $0\leq \a <1$ every Riemannian metric of
the following family
$$
\begin{array}{c}
\bar{g}_\alpha= \frac{r^4 (r^2-\alpha^2) (r^2 +\alpha^2)}{r^8-2
\alpha^4 (r^4-1) -1} dr^2 + \frac{r^8-2 \alpha^4 (r^4-1) -1}{r^2
(r^2-\alpha^2) (r^2 +\alpha^2)} \eta_1^2 + r^2 (\eta_2^2+
\eta_3^2)
\\ \\
+ (r^2+\alpha^2) (\eta_4^2+\eta_5^2) + (r^2-\alpha^2)
(\eta_6^2+\eta_7^2),
\end{array}\eqno{(**)}
$$
is smooth Riemannian metric on $H/\Z_2$ with $SU(4)$-holonomy. For
$\a=0$ metric $(**)$ is isometric to the Calabi metric \cite{Calabi}
with $SU(4)$-holonomy; for $\alpha=1$ metric $(**)$ is isometric to
the Calabi metric \cite{Calabi} with holonomy group $Sp(2) \subset
SU(4)$ on $T^*\C P^2$. }

\vskip0.2cm

Notice that the form of metric $\bar{g}_\alpha$ in Theorem 1 for
$a=0$ and $\a=1$ differs from \cite{Calabi}; the Calabi metrics of
such form were described in \cite{Page-Pope} and \cite{CGLP}.

This result was obtained while systematic studying the metrics
$(*)$ with holonomy $Spin(7)$ using technic introduced in
\cite{Baz1} and exploited then in \cite{Baz-Malk,Baz2}: metric
$(*)$ can be constructed using arbitrary 7-dimensional 3-Sasakian
manifold $M$ and will possesses natural $Spin(7)$-structure. The
parallelness of this structure can be reduced to the system of
nonlinear ODE's:
$$
\begin{array}{l}
A_1'=\frac{(A_2-A_3)^2-A_1^2}{A_2A_3} +
\frac{A_1^2(B^2+C^2)}{B^2C^2},\\
A_2'=\frac{A_1^2-A_2^2+A_3^2}{A_1A_3}-\frac{B^2+C^2-2A_2^2}{BC},
\\
A_3'=\frac{A_1^2+A_2^2-A_3^2}{A_1A_2}-\frac{B^2+C^2-2A_3^2}{BC},
\\
B'=-\frac{CA_1+BA_2+BA_3}{BC} -
\frac{(C^2-B^2)(A_2+A_3)}{2A_2A_3C}, \\
C'=-\frac{BA_1+CA_2+CA_3}{BC} -
\frac{(B^2-C^2)(A_2+A_3)}{2A_2A_3B}.
\end{array}\eqno{(***)}
$$
Notice that the system $(***)$ in the case of $B=C$ was completely
explored in \cite{Baz1,Baz2}. To obtain smooth metric $(*)$ one
should resolve the conic singularity on $\bar{M}$ in one of two
possible ways give rise to a spaces ${\cal M}_1$ or ${\cal M}_2$.
This scheme is described in the second section of our article.
Solution $(**)$ on ${\cal M}_2/ \Z_2$ can be constructed by
integrating system $(***)$ assuming $A_2=-A_3$ (section three).
Notice that metrics presented in Theorem 1 have cohomogeneity one
with respect to natural action of Lie group $SU(3)$. We also should
remark that cohomogeneity one spaces with parallel
$Spin(7)$-structures were studied in \cite{Reid}.

The proof of the following theorem is rather complicated and will be
published separately.

\vskip0.2cm

{\bf Theorem 2.} {\it Let $M$ be a $7$-dimensional compact
$3$-Sasakian manifold with K\"{a}hler ${\cal O}$ and let $p=2$ if a
general fibre of $3$-Sasakian bundle equals to $SO(3)$ or $p=4$ if
the general fibre is $SU(2)$. Then there exist following complete
regular Riemannian metrics $\bar{g}$ of the form $(*)$ with holonomy
group $H \subset Spin(7)$ on orbifold ${\cal M}_2/\Z_p$:

1) if $A_1(0)=0$, $-A_2(0)=A_3(0)>0$ and $2A_2^2(0)=B^2(0)+C^2(0)$
then metric $\bar{g}$ of the form $(*)$ possesses holonomy $SU(4)
\subset Spin(7)$ and homothetic to the one of the family $(**)$;

2) if $A_1(0)=0$ and $-A_2(0)=A_3(0)<B(0)=C(0)$ then there exist
regular ALC-metric $\bar{g}$ of the form $(*)$ with holonomy
$Spin(7)$ found in \cite{Baz1}. At infinity this metrics tend to a
product of cone over a twistor space ${\cal Z}$ and a circle.

Moreover, any complete metric $(*)$ on ${\cal M}_2/\Z_p$ with
described $Spin(7)$-structure and with holonomy $H \subset Spin(7)$
is isometric to the one of the above-mentioned. }

\vskip0.2cm

\section[]{Description of the $Spin(7)$-structure on a cone over a $3$-Sasakian manifold.}

In this section we will describe, as briefly as possible, the spaces
on which we construct metrics with parallel $Spin(7)$-structure. In
our further notations and definitions concerned with $3$-Sasakian
manifolds we will follow the paper \cite{Baz1}. For more details on
$3$-Sasakian manifolds look paper \cite{Boyer-Galicki}.

Let $\bar{M} = \mathbb{R}_+ \times M$, $t \in \R_+ = (0, \infty)$ be
the cone over smooth closed Riemannian manifold $(M,g)$ with the
metric $dt^2+t^2g$. Manifold $M^m$ is called $3$-Sasakian if the
cone metric is hyperk\"{a}hler, that is holonomy group of $\bar{M}$
is contained in $Sp(\frac{m+1}{4})$. Then there are three parallel
complex structures $J^1,J^2,J^3$ on $\bar{M}$ satisfying
$J^jJ^i=-\delta^{ij}+ \varepsilon_{ijk}J^k$. Identifying $M$ with a
``base'' of the cone $\bar{M} \cap \{t=1\}$ consider vector fields
on $M$: $\xi^i=J^i(\partial_t)$, $i=1, 2, 3$. Fields $\xi_i$ are
called characteristic fields of the $3$-Sasakian manifold $M$, and
dual $1$-forms $\eta_i$ are called characteristic forms. It can be
shown that fields $\xi_i$ generate the Lie algebra ${\bf su}(2)$
with respect to the Lie bracket as an operation on vector fields.
Hence there is a bundle $\pi: M \rightarrow {\cal O}$ with
$SU(2)=S^3$ or $SO(3)=\R P^3$ as a general fibre over certain
$4$-dimensional quaternionic K\"{a}hler orbifold ${\cal O}$. Denote
a bundle of horizontal (with respect to $\pi$) vector fields on $M$
as ${\cal H}$.

Define following $2$-forms on $M$:
$$
\omega_i=d\eta_i+\sum_{j,k}\varepsilon_{ijk}\eta_j \wedge \eta_k,
\ i=1,2,3.
$$
One can check \cite{Baz1} that forms $\omega_i$ generate a subspace
in $\Lambda^2 \mathcal{H^*}$, so we can fix orthonormal system of
$1$-forms $\eta_4,\eta_5,\eta_6,\eta_7$ on the $\mathcal{H}$ such
that
$$
\begin{array}{l}
\omega_1=2( \eta_4 \wedge \eta_5 - \eta_6 \wedge \eta_7), \\
\omega_2=2( \eta_4 \wedge \eta_6 - \eta_7 \wedge \eta_5), \\
\omega_3=2( \eta_4 \wedge \eta_7 - \eta_5 \wedge \eta_6).
\end{array}
$$
Let $\R^8$ be a standard euclidean space with coordinates $x^0,
\ldots x^7$. Denote $e^{ijkl}= dx^i \wedge dx^j \wedge dx^k \wedge
dx^l$ and define self-adjoint $4$-form on $\R^8$ as follows:
$$
\begin{array}{l}
\Phi_0=e^{0123}+e^{4567}+e^{0145}-e^{2345}-e^{0167}+e^{2367}+e^{0246}+
\\
e^{1346}-e^{0275}+e^{1357}+e^{0347}-e^{1247}-e^{0356}+e^{1256}.
\end{array}
$$
It is known that the group of linear transformation of $\R^8$
preserving $\Phi_0$ is isomorphic to $Spin(7)$. Moreover, in
addition group $Spin(7)$ preserve orientation and metric
$g_0=\sum_{i=0}^7 (e^i)^2$. Let $N$ be an oriented Riemannian
$8$-manifold. A differential form $\Phi \in \Lambda^4 N$ generates
$Spin(7)$-structure on $N$ if there is an orientation-preserving
isometry $\phi_p : T_p N \rightarrow \R^8$ such that
$\phi_p^*\Phi_0=\Phi|_p$ for any $p \in N$. If form $\Phi$ is
parallel then holonomy of Riemannian manifold $N$ is reduced to a
subgroup $Spin(7) \subset SO(8)$, i.e. $Hol(N) \subset Spin(7)$.
It's known \cite{Gray1} that the parallelness of $\Phi$ is
equivalent to it's closeness:
$$
d \Phi=0.
$$
Consider $Spin(7)$-structure on the $\bar{M}$ defined by the
following form
$$\begin{array}{l}
\Phi=e^{0123}+C^2B^2\eta_4\wedge\eta_5\wedge\eta_6\wedge\eta_7 +
\frac{B^2+C^2}{4}(e^{01}-e^{23})\wedge \omega_1 \\
+\frac{B^2-C^2}{4}(e^{01}-e^{23})\wedge \omega
+\frac{BC}{2}(e^{02}-e^{31})\wedge \omega_2
+\frac{BC}{2}(e^{03}-e^{12})\wedge \omega_3,
\end{array}$$
where
$$
\begin{array}{l}
e^0=dt, \\
e^i=A_i \eta_i, i=1,2,3, \\
e^j=B \eta_j, j=4,5, \\
e^k=C \eta_k, k=6,7,
\end{array}
$$
$A_1(t), A_2(t), A_3(t), B(t), C(t)$ be some smooth functions. It's
easy to understand that form $\Phi$ corresponds to metric
$$
dt^2+ A_1(t)^2 \eta_1^2 + A_2(t)^2 \eta_2^2 + A_3(t)^2 \eta_3^2+
B(t)^2 (\eta_4^2+\eta_5^2) + C(t)^2 (\eta_6^2+\eta_7^2) \eqno{(1)}
$$
on $\bar{M}$. We will suppose that quaternionic K\"{a}hler orbifold
${\cal O}$ is K\"{a}hler and co-frame $\eta_i$, $i=4,5,6,7$ can be
chosen such that form $\omega=2( \eta_4 \wedge \eta_5 + \eta_6
\wedge \eta_7)$ generates K\"{a}hler structure on ${\cal O}$ and, in
particularly, is closed. This assumption encloses exterior algebra
of the forms under the consideration and allows to derive a correct
system of ODE's on the functions $A_i, B, C$. Notice that without
K\"{a}hlerness assumption of the ${\cal O}$ it is necessary to put
$B=C$ to enclose the forms algebra. Using expressions on exterior
algebra from \cite{Baz1}
$$
\begin{array}{l}
d e^0=0, \\
d e^i=\frac{A_i'}{A_i} e^0 \wedge e^i+ A_i \omega_i- \frac{2
A_i}{A_{i+1} A_{i+2}} e^{i+1} \wedge e^{i+2}, \ i=1,2,3 \mbox{\ mod \ } 3, \\
d \omega_i = \frac{2}{A_{i+2}} \omega_{i+1} \wedge e^{i+2}
-\frac{2}{A_{i+1}} e^{i+1} \wedge \omega_{i+2},  \ i=1,2,3 \mbox{\
mod \ } 3,
\end{array}
$$
and condition $d \omega =0$, after some computations we immediately
obtain

\vskip0.2cm

{\bf Lemma 1.}  {\it The parallelness of the $\Phi$ is equivalent to
the following system of ODE's:
$$
\begin{array}{l}
A_1'=\frac{(A_2-A_3)^2-A_1^2}{A_2A_3} +
\frac{A_1^2(B^2+C^2)}{B^2C^2},\\
A_2'=\frac{A_1^2-A_2^2+A_3^2}{A_1A_3}-\frac{B^2+C^2-2A_2^2}{BC},
\\
A_3'=\frac{A_1^2+A_2^2-A_3^2}{A_1A_2}-\frac{B^2+C^2-2A_3^2}{BC},
\\
B'=-\frac{CA_1+BA_2+BA_3}{BC} -
\frac{(C^2-B^2)(A_2+A_3)}{2A_2A_3C}, \\
C'=-\frac{BA_1+CA_2+CA_3}{BC} -
\frac{(B^2-C^2)(A_2+A_3)}{2A_2A_3B}.
\end{array}\eqno{(2)}
$$ }

\vskip0.2cm

For $B=C$ we obtain a system explored in \cite{Baz1}:
$$
\begin{array}{l}
A_1'=\frac{2 A_1^2}{B^2}+\frac{(A_2-A_3)^2-A_1^2}{A_2A_3},\\
A_2'=\frac{2 A_2^2}{B^2}+\frac{(A_3-A_1)^2-A_2^2}{A_1A_3},\\
A_3'=\frac{2 A_3^2}{B^2}+\frac{(A_1-A_2)^2-A_3^2}{A_1A_2},\\
B'=-\frac{A_1+A_2+A_3}{B}.
\end{array}\eqno{(2')}
$$
To construct a smooth Riemannian metric on manifold (or orbifold) it
is necessary to set initial data for $(2)$. It was shown in
\cite{Baz1} that there are two possible ways to resolve the conic
singularity of the $\bar{M}$, corresponding to two spaces ${\cal
M}_1$ and ${\cal M}_2$. Further we describe the space ${\cal M}_2$
on which we will search for the metric with holonomy $H \subset
Spin(7)$.

Let $S\simeq S^1$ --- subgroup in $SU(2)$ (or in $SO(3)$)
integrating one of the characteristic Killing fields, for instance
$\xi_1$. Hence there is principal bundle $\pi': M \rightarrow {\cal
Z}$ with structure group $S$, where ${\cal Z}= M/S$ --- twistor
space. Consider natural action $S$ on $\R^2=\C$: $ e^{i\phi} \in S:
z \rightarrow z e^{i \phi}$ and associate with $\pi'$ fibred space
${\cal M}$ with fibre $\C$ with respect to considered action. Thus,
orbifold ${\cal Z}$ is imbedded into ${\cal M}_2$ as a zero section,
and ${\cal M}_2 \backslash {\cal Z}$ is foliated on the spherical
subbundles diffeomorphic to the $M$ and collapsing to the zero
section ${\cal Z}$ for $t \rightarrow 0$.

Now let $p \in \N$, $Z_p \subset S$. Group $Z_p$ acts on ${\cal
M}_2$ by isometries. Consequently, orbifold  ${\cal M}_2/ \Z_p$ is
well-defined and it is manifold if and only if ${\cal M}_2$ is
manifold. It easy to understand that ${\cal M}_2 / \Z_p$ is bundle
with fibre $\C$ associated with principal bundle $\pi':M
\rightarrow {\cal Z}$ with respect to the action $e^{i\phi} \in S:
z \rightarrow z e^{i p \phi}$.

Notice that if $3$-Sasakian manifold is regular (i.e. fibration on
$3$-dimensional $3$-Sasakian fibres is regular) then all fibres of
the $\pi$ are isometric to the $S^3=SU(2)$ (or $SO(3)$) and
orbifolds ${\cal O}$, ${\cal Z}$ and ${\cal M}_2$ turn out to be
smooth manifolds. It's known \cite{Boyer-Galicki} that it is
possible only if $M$ is isometric to $S^7$, $\R P^7$ or
$N_{1,1}=SU(3)/T_{(1,1)}$. But only the Aloff-Wallach space
$N_{1,1}$ has K\"{a}hler base ${\cal O}$ and new metrics can be
constructed only in this case.

The following lemma proposes necessary conditions on the functions
$A_i, B, C$ for the solution of $(2)$ to determine a smooth metric
$(1)$ on the ${\cal M}_2$.

\vskip0.2cm

{\bf Lemma 2.} {\it Let $(A_1(t),A_2(t),A_3(t),B(t),C(t))$ is a
$C^{\infty}$ solution of the system $(2)$, $t \in [0, \infty)$. Let
$p=4$ for general fibre to be isometric to $SO(3)$, and $p=2$ for
$SU(2)$. Then the metric $(1)$ can be continued to a smooth metric
on the ${\cal M}_2/ \Z_p$ if and only if the following conditions
hold

(1) $A_1(0)=0, |A_1'(0)|=4$;

(2) $A_2(0)=-A_3(0) \neq 0, A_2'(0)=A_3'(0)$,

(3) $B(0) \neq 0, B'(0)=0$;

(4) $C(0) \neq 0, C'(0)=0$;

(5) the functions $A_1,A_2,A_3,B,C$ have fixed sign on the
interval $(0,\infty)$. }

\vskip0.2cm

Lemma 2 was proven in \cite{Baz1} for the case $B=C$. The proof can
be generalized without any serious modifications, only the following
should be noticed: in \cite{Baz1} choice of the field $\xi_i$
corresponding to collapsing of the $S$ have no importance for
constructing ${\cal M}_2$ because of additional symmetries of the
system $(2')$. But some simple calculations show that in general
case it is necessary to choose $\xi_1$ as a generator for $S$,
because the only function $A_1$ can vanish at the initial point
$t=0$.

\section[]{Construction of the explicit solutions on $\mathcal{M}_2$.}

In system $(2)$ we make a substitution $A_2=-A_3$. Then summing the
second and the third equations we obtain $B^2+C^2=2 A_2^2$ and
subtracting the fourth and the fifth equations we get
$(B^2-C^2)'=0$. So we can assume that $B^2=A_2^2+\alpha^2$,
$C^2=A_2^2-\alpha^2$ for some constant $\alpha \geq 0$ and system
$(2)$ is reduced to
$$
\begin{array}{c}
A_1'=-4+\frac{A_1^2}{A_2^2} + 2 \frac{A_1^2 A_2^2}{A_2^4-\alpha^4},
\\
(A_2^2)'=-A_1.
\end{array}
$$

This system can be integrated easily. To be more precise, consider a
new variable $\rho$ such that $d\rho=-2A_1 dt$. Adding some constant
to $\rho$ we can assume that $A_2^2=\rho$ and putting $A_1^2=F$ we
obtain
$$
\frac{dF}{d\rho} + F G = 4,
$$
where
$$
G(\rho)= \frac{1}{\rho}+\frac{1}{\rho-\alpha^2} +
\frac{1}{\rho+\alpha^2}.
$$
That system is easy to solve (using integrating multiplier) and
putting $r^2=\rho$ we eventually obtain
$$
 F=\frac{r^8-2 \alpha^4 r^4
+ \beta}{r^2 (r^4-\alpha^4)},
$$
where $\beta$ is an integrating constant. Consequently, metric $(1)$
becomes
$$
\begin{array}{c}
\bar{g}= \frac{r^4 (r^2-\alpha^2) (r^2 +\alpha^2)}{r^8-2 \alpha^4
r^4 +\beta} dr^2 + \frac{r^8-2 \alpha^4 r^4 +\beta}{r^2
(r^2-\alpha^2) (r^2 +\alpha^2)} \eta_1^2 + r^2 (\eta_2^2+ \eta_3^2)
\\ \\
+ (r^2+\alpha^2) (\eta_4^2+\eta_5^2) + (r^2-\alpha^2)
(\eta_6^2+\eta_7^2).
\end{array}
$$
To obtain a regular metric on $\mathcal{M}_2/ \Z_p$ it is necessary
for polynomial $r^8-2 \alpha^4 r^4 + \beta$ to have roots, moreover,
the maximal root $r_0$ must be greater than $\alpha$. In this case
metric is well-defined for $r \geq r_0$. It is obvious that
replacing the starting metric by a homothetic metric we can
normalize the maximal root $r_0$ as $r_0=1$. Some trivial
calculations show that this means $0 \leq \alpha <1$ and
$\beta=2\alpha^4-1$. Thus, metric $(1)$ becomes
$$
\begin{array}{c}
\bar{g}_\alpha= \frac{r^4 (r^2-\alpha^2) (r^2 +\alpha^2)}{r^8-2
\alpha^4 (r^4-1) -1} dr^2 + \frac{r^8-2 \alpha^4 (r^4-1) -1}{r^2
(r^2-\alpha^2) (r^2 +\alpha^2)} \eta_1^2 + r^2 (\eta_2^2+ \eta_3^2)
\\ \\
+ (r^2+\alpha^2) (\eta_4^2+\eta_5^2) + (r^2-\alpha^2)
(\eta_6^2+\eta_7^2),
\end{array}\eqno{(3)}
$$
where $0 \leq \alpha <1$, $r\geq 1$. Immediate checking of the
conditions of Lemma $(2)$ shows that $(3)$ is a family of a smooth
metrics on $\mathcal{M}_2 /\Z_p$ for $r \geq 1$ and $0 \leq \alpha
<1$. Moreover, $\bar{g}_0$ coincides with Calabi's metric with
holonomy $SU(4)$ constructed in \cite{Calabi}.

It follows from Lemma $(1)$ that holonomy $Hol (\bar{g}_\alpha)$ of
the metric $(3)$ is contained in $Spin(7)$. Let us consider the
following $2$-form:
$$
\bar{\Omega}_1=-e^0\wedge e^1 +e^2\wedge e^3 +e^4\wedge e^5
-e^6\wedge e^7.
$$
Obviously, this form is agreed with metric $(1)$ and simple
calculation shows that it is closed exactly for $A_2=-A_3$; it means
that this form is K\"{a}hler form of the metric $(3)$. Hence
$Hol(\bar{g}_\alpha) \subset SU(4)$.

Further we explore more carefully the case when the metrics
$\bar{g}_\alpha$ are defined on a smooth manifold, i.e. when
$M=N_{1,1}$. Firstly we make some useful notation for the following
subgroups in $SU(3)$:
$$
S_{1,1}=\{ \diag(z,z,\bar{z}^{2}) | z \in S^1\subset \C\},
$$
$$
T=\{ \diag(z_1,z_2,\bar{z}_1\bar{z}_2) | z_1,z_2 \in S^1\subset
\C\},
$$
$$
K_1= \left\{ \left( \begin{array}{cc} \det (\bar{A}) & 0 \\ 0 & A
\end{array} \right) | A \in U(2) \right\},
$$
$$
K_3= \left\{ \left( \begin{array}{cc} A & 0 \\ 0 & \det (\bar{A})
\end{array} \right) | A \in U(2) \right\},
$$
Let $\C^3$ be a 3-dimensional complex space, let $S^5 \subset \C^3$
be a unit sphere and let unit circle $S^1$ acts diagonally on $\C^3$
and it's associated spaces. This action defines an equivalence class
noted further by square brackets: $[u,v]$, $[u]$ etc.

Let $\tilde{E} =\{ (u_1, u_2) | |u_1|=1, \langle u_1, u_2 \rangle_\C
=0 \} \subset S^5 \times \C^3$. Consider a diagonal action of a
circle $S^1$ on the $\tilde{E}$ and a projection $\tilde{\pi}_1:
(u_1,u_2) \mapsto u_1$ on the $\tilde{E}$ to $S^5$ which is a bundle
with fibre $\C^2$. The spherical fibration of the $\tilde{\pi}_1$ is
$\tilde{E}^1=\{ (u_1, u_2) \in E | |u_1|=|u_2|=1, \langle u_1, u_2
\rangle_\C=0 \}$ diffeomorphic to $SU(3)$. The bundle $\tilde{\pi}$
(with $S^1$-action) induces a vector bundle $\pi_1: E=\tilde{E}/S^1
\rightarrow \C P^2$ with fibre $\C^2$ and a spherical fibration
$E^1=\tilde{E}^1/S^1=SU(3)/S_{1,1}=N_{1,1} \rightarrow \C
P^2=SU(3)/K_1$. It is easy to understand that $\pi_1$ can be
identified with cotangent space $T^* \C P^2 \rightarrow \C P^2$.

Similarly, we consider space $\tilde{H}= \{ (u_1, u_2, [u_3]) |
|u_1|=|u_3|=1, \langle u_i, u_j \rangle_\C = 0, i,j=1, 2, 3 \}
\subset S^5 \times \C^3 \times \C P^2$ and projection
$\tilde{\pi}_2: (u_1,u_2,[u_3]) \mapsto (u_1, [u_3])$ of the space
$\tilde{H}$ to the space $\tilde{F}= \{ (u_1, [u_3]) |
|u_1|=|u_3|=1, \langle u_1, u_3 \rangle_\C =0 \}$ with fibre $\C$.
The spherical fibration of the $\tilde{\pi}_2$ coincides with
$\tilde{H}^1=\{ (u_1, u_2, [u_3]) | \langle u_i, u_j \rangle_\C=0,
|u_i|=1, i,j=1,2, 3 \}$ and obviously can be identified with
$SU(3)=\tilde{E}^1$. The bundle $\tilde{\pi}_2$ (with the same
$S^1$-action) induces a bundle $\pi_2: H=\tilde{H}/S^1 \rightarrow
F=\tilde{F}/S^1$ with fibre $\C$ which's spherical fibration
coincides with $E^1=\tilde{E}^1=N_{1,1} \rightarrow SU(3)/T$. The
base of the $\pi_2$ is complex flag-manifold $F=SU(3)/T$ which can
be regarded as follows:
$$
F=\{ ([u_1],[u_3])| u_i \in \C^3, |u_i|=1, \langle u_1, u_3
\rangle_\C=0, i=1,3 \}.
$$
The complex line bundle $\pi_2: H \rightarrow F$ is said to be {\it
canonical} bundle over the manifold $F$ of complex flags in $\C^3$.

Thus canonical bundle over $F$ and cotangent bundle over $\C P^2$
have common spherical subbundle $N_{1,1}$, which fibered in two
different ways. It is well-known fact \cite{Boyer-Galicki}, that the
Aloff-Wallach space $M=N_{1,1}$ admits $3$-Sasakian structure with a
twistor fibration coinciding with the bundle $\pi_2: N_{1,1}
\rightarrow F={\cal Z}$, and with $3$-Sasakian fibration coinciding
with $SO(3)$-fiber projection $\pi_2':N_{1,1} \rightarrow
SU(3)/K_3=\C P_2={\cal O}$. It is evidently that in this case ${\cal
M}_2$ coincide with above considered fibred space $H$ over flag
manifold $F$. If $0 \leq \alpha <1$ then (3) is smooth complete well
defined metric on $H/\Z_2$, that is on the space of complex line
bundle $\pi_2 \otimes \pi_2$. Under the assumption $\alpha=1$ metric
(3) is reduced to the Calabi metric \cite{Calabi} on $E=T^* \C P^2$.

\vskip0.2cm

{\bf Theorem 1.} {\it Let $M=N_{1,1}$ be the Aloff-Wallach space,
then explicit Riemannian metrics $\bar{g}_\alpha$ described in (3)
are pairwise nonhomothetic smooth complete metrics with the
following properties:

- if $0 \leq \alpha<1$, then metric $\bar{g}_\alpha$ is smooth
metric on space $H/\Z_2$ of tensor square of canonical bundle
$\pi_2:H\rightarrow F$ over the flag manifold $F$ in $\C^3$ and has
holonomy group $SU(4)$. Metric $\bar{g}_0$ coincides with the Calabi
metric \cite{Calabi};

- metric $\bar{g}_1$ has holonomy group $Sp(2) \subset SU(4)$ and
coincides with the hyperk\"ahler Calabi metric \cite{Calabi} on $T^*
\C P^2$.

}

\vskip0.2cm

{\bf Proof.} To see hyperk\"ahlerness of $\bar{g}_1$, it is
sufficiently to consider additional pair of K\"ahler forms,
generating together with $\bar{\Omega}_1$ a hyperk\"ahler structure:
$$
\bar{\Omega}_2=e^0\wedge e^2 +e^1\wedge e^3 - e^4\wedge e^6
+e^7\wedge e^5 = e^0\wedge e^2 +e^1\wedge e^3 -
\frac{BC}{2}\omega_2,
$$
$$
\bar{\Omega}_3=-e^0\wedge e^3 +e^1\wedge e^2 - e^4\wedge e^7
+e^5\wedge e^6 = -e^0\wedge e^3 +e^1\wedge e^2 -
\frac{BC}{2}\omega_3.
$$
Immediate calculations show that forms $\bar{\Omega}_2$,
$\bar{\Omega}_2$ are close exactly in case of $\alpha=1$. This gives
reduction of holonomy group to $Sp(2) \subset SU(4)$ for metric
$\bar{g}_1$.

For completing the proof we need to show that $\bar{g}_\alpha$ is
not hyperk\"ahler in case $0 \leq \alpha<1$. Indeed, if $Hol
(\bar{g}_\alpha) = Hol ({\cal M}_2/\Z_2) \subset Sp(2)$ for $0 \leq
\alpha < 1$, then limit metric has the same property: $Hol (\bar{M}/
\Z_2) \subset Sp(2)$. However, it is clear that after factorizing
cone $\bar{M}$ by group $\Z_2$, the generator of $\Z_2$ has to be
added to holonomy group of $\bar{M}$. This generator corresponds to
transformation $\H^2 \rightarrow \H^2: (q_1, q_2) \mapsto
(q_1',q_2')$, where $q_l=u_l+v_l j$ and $q_l'=-u_l +v_l j$, $u_l,v_l
\in \C$, $l=1,2$. Above transformation is an element of $SU(4)$ but
clearly does not belong to $Sp(2)$. Consequently $Hol(\bar{M}/\Z_2)$
is not a subgroup of $Sp(2)$ and $Hol(\bar{g}_\alpha)=SU(4)$.

\section[]{Some generalizations.}

The family of metrics $(**)$ could be generalized. It turns out that
in every dimension $4n, n\geq2$ there is a continuous family of
metrics $G_\alpha, \alpha \in [0,1]$ "connecting" the Calabi metric
$G_0$ with $SU(2n)$-holonomy and the Calabi metric $G_1$ with
$Sp(n)$-holonomy \cite{Malkovich}.

\end{document}